\documentclass[12pt]{article}
\usepackage{amssymb, amsthm, amsmath,eucal,
pstricks,pst-plot,pst-node,eucal}
\usepackage[dvips]{graphics}

\newpsobject{showgrid}{psgrid}{subgriddiv=1,griddots=5,gridlabels=0pt}

\newcommand{\Q}{\mathbb{Q}}
\newcommand{\C}{\mathbb{C}}
\newcommand{\T}{\mathbb{T}}

\newcommand{\Zt}{\mathbb{Z}/2}

\newcommand{\cO}{{\mathcal{O}}}
\newcommand{\cC}{{\mathcal{C}}}
\newcommand{\cB}{{\mathcal{B}}}

\newcommand{\cBp}{\cB^\times}
\newcommand{\fb}{\mathbf{f}}
\newcommand{\pb}{\mathbf{p}}

\newcommand{\sA}{{\mathsf{A}}}
\newcommand{\sS}{{\mathsf{S}}}
\newcommand{\sB}{{\mathsf{B}}}
\newcommand{\sC}{{\mathsf{C}}}
\newcommand{\sG}{{\mathsf{G}}}

\newcommand{\hMon}{\widehat{\Mon}}

\newcommand{\gl}{{\mathfrak{gl}}}
\newcommand{\fq}{{\mathfrak{q}}}

\newcommand{\ran}{\right\rangle}
\newcommand{\lan}{\left\langle}

\DeclareMathOperator{\Aut}{Aut}

\DeclareMathOperator{\wnd}{wind}

\DeclareMathOperator{\Mon}{\mathbf{M}}
\DeclareMathOperator{\Lift}{Lift}
\DeclareMathOperator{\wt}{wt}
\DeclareMathOperator{\Cliff}{Cliff}

\numberwithin{equation}{section}

\newtheorem{Theorem}{Theorem}
\newtheorem{Lemma}{Lemma}
\newtheorem{Corollary}[Lemma]{Corollary}
\newtheorem{Proposition}[Lemma]{Proposition}

\begin{document}
\title{The theta characteristic of a branched covering}
\author{Alex Eskin, Andrei Okounkov, and Rahul Pandharipande}
\date{August 2003}
\maketitle

\begin{abstract}
We give a group-theoretic
description of the parity 
of a pull-back of a theta
characteristic under a 
branched covering. It 
involves lifting 
 monodromy of the covering 
to the semidirect product 
of the symmetric and Clifford
groups, known as the Sergeev 
group. As an application, we
enumerate torus coverings
with respect to their 
ramification and parity and, 
in particular, 
show that the corresponding 
all-degree generating 
functions are quasimodular 
forms. 
\end{abstract}

\section{Introduction}

\subsection{Statement of the problem}

\subsubsection{}

Let
\begin{equation}
  \label{deff}
  f: \cC \to \cB
\end{equation}
be a degree $d$ map between smooth algebraic curves. 
Let $L$ be a theta characteristic of $\cB$, that is,
a divisor such that 
$$
2 L = K_\cB\,,
$$
where $K_\cB$ is the canonical class of $\cB$. For example,
if $\cB$ is a rational curve then any divisor of degree
$-1$ is a theta characteristic. We will allow disconnected
curves $\cC$. 

Throughout
this paper we will assume that $f$ has only 
\emph{odd ramifications},
that is, its local form near any point of $\cC$ is 
$$
z\mapsto z^k\,,\quad k=1,3,5,\dots\,.
$$
Under this assumption, one can define the pull-back 
of $L$ under $f$ as follows. Let $\omega$ be a meromorphic
$1$-form on $\cB$ such that $(\omega)=2L$. All zeros of 
the $1$-form $f^*\omega$ have even multiplicity and so 
we can set
$$
f^*L = \frac12 \, \left(f^*\omega\right) \,.
$$

\subsubsection{}

A theta characteristic $L$ has 
a deformation invariant, namely its 
parity
\begin{equation}
  \label{defp}
  p(L) = \dim  \cO(L)\!\!\!\!\mod 2 \,.
\end{equation}
This invariant 
has been  studied
from various points of view both classically and 
recently, see \cite{At,Har,Mum}.

The main question addressed in this 
paper is  how the parity of $f^*L$ 
is related to the parity of $L$ and 
the topological data of the branched
covering $f$. The answer is given 
in Theorem \ref{t1}. It involves the 
combinatorics of the Sergeev group $\sC(d)$
defined in Section \ref{sSe}.

\subsubsection{}

The particular case of 
a special importance to us will be the
case when $\cB=\T^2$ is a standard torus
and $L=0$. Using the representation 
theory of the Sergeev group, we enumerate
the branched coverings 
\begin{equation}
  \label{ftoT}
  f: \cC \to \T^2\,,
\end{equation}
with given ramification 
with respect to the parity of their theta
characteristic. This is done in Theorem 
\ref{t2}. 

As a corollary, we prove that 
the natural all-degree generating
function for these numbers is a 
quasi-modular form, refining a result
of \cite{BO}, see Corollary \ref{c1}.  This quasimodularity 
is a very useful property for 
the $d\to\infty$ asymptotic enumeration
of coverings along the lines of \cite{EO}. 
It gives a way to compute the volumes
of the corresponding strata in the moduli
space of holomorphic differentials \cite{KZ}.
This application was the main motivation 
for the present work.

\subsection{Johnson-Thurston formula}

\subsubsection{}

Every theta characteristic $L$ determines a 
quadratic form $q_L$ on the group 
$$
J_2(\cB) \cong H^1(\cB,\Zt)
$$
of elements of order 2 in the Jacobian of $\cB$ 
by the following rule:
\begin{equation}
q_L(\gamma) = p(L\otimes \gamma) + p(L)
\label{defq}
\end{equation}
The Arf invariant of this quadratic form is $p(L)$ itself,
meaning that 
\begin{equation}
  \label{arf}
  (-1)^{p(L)} = 2^{-g(\cB)} \, \sum_{\gamma\in J_2(\cB)} 
(-1)^{q_L(\gamma)} \,,
\end{equation}
assuming $\cB$ is connected. Here $g(\cB)$ is the genus
of $\cB$. 

Throughout the paper, we identify $H^1(\cB)$ and $H_1(\cB)$
using the intersection pairing. 
By the Riemann-Mumford 
relation \cite{ACGH}, the bilinear form associated 
with the quadratic form $q_L$ is the intersection 
form.

\subsubsection{}

We will use a different, purely topological, description 
of the form $q_L$ which is due to Johnson \cite{Jo} and 
Thurston \cite{ACGH}.  Let $\omega$ be a meromorphic $1$-form 
such that $(\omega)=2L$. Away from the support of $(\omega)$,
we can introduce a local coordinate $z$ such that 
$$
\omega = dz \,.
$$
The transition functions between  different charts will then
be translations of $\C$, thus making $\cB$ a 
so-called \emph{translation surface}. 

Let $\Gamma$ be an immersed closed curve on $\cB$ that avoids the
support of $(\omega)$. The translation surface structure
determines the number of full rotations of the tangent vector 
to $\Gamma$. Taken modulo 2 this number is independent of 
the orientation of $\Gamma$. We will denote it by
$$
\wnd(\Gamma)\in \Zt \,.
$$ 
The Johnson-Thurston formula can now be stated as follows
\begin{align}
  \label{TJ}
  q_L(\gamma) = \wnd(\Gamma) \,\,+ \,\, 
&\textup{\# of components of $\Gamma$} \\
+ \,\, &\textup{\# of self-intersections of $\Gamma$} \notag
\,,
\end{align}
where $\Gamma$ is a representative of the
class $\gamma$. Note, in particular, that the right-hand
side of \eqref{TJ} depends only on the homology class mod $2$ of 
the curve $\Gamma$, which is elementary to check directly.

\subsubsection{}

Given a map as in \eqref{deff} and a class $\gamma\in H_1(\cC,\Zt)$,
we set
\begin{equation}
  \label{defDq}
  \Delta_f (\gamma) = q_{f^*L}(\gamma) - q_L(f_* \gamma) \,.
\end{equation}
It is clear from \eqref{TJ} that this does not depend on the 
choice of the theta-characteristic $L$ and simply equals
the number of ``new'' self-intersections of $f_*\gamma$, that
is, 
\begin{align}
  \label{newS}
  \Delta_f (\gamma) = & \textup{\# of self-intersections of $f(\Gamma)$} - \\
& \textup{\# of self-intersections of $\Gamma$}\notag\,,
\end{align}
where $\Gamma$ is a representative of $\gamma$. 

\subsection{The Sergeev group}

\subsubsection{}

The quadratic form \eqref{newS} is a purely combinatorial
invariant of a branched covering \eqref{deff} and a 
homology class $\gamma\in H_1(\cC,\Zt)$. Our first result 
is a group-theoretic interpretation of this invariant.

Every element $\gamma\in J_2(\cC)$ defines a unramified
double cover $\cC_\gamma \to \cC$, and hence, by composing with
$f$,  a degree $2d$ 
cover 
\begin{equation}
  \label{fga}
f_\gamma: \cC_\gamma \to \cB \,.
\end{equation}
Let 
\begin{equation}
  \label{defsi}
  \sigma: \cC_\gamma \to \cC_\gamma
\end{equation}
be the involution permuting the sheets of $\cC_\gamma \to \cC$. 
The monodromy of the branched covering 
\eqref{fga} lies in the centralizer of a fixed-point-free
involution $\sigma$ inside the symmetric group $\sS(2d)$. 

In other words, the monodromy of $f_\gamma$ defines a 
homomorphism 
\begin{equation}
\Mon_{f_\gamma}: \pi_1(\cBp) \to  \sB(d) \,,
\label{phig}
\end{equation}
where
\begin{equation}
  \label{cBp}
  \cBp = \cB \setminus \{\textup{branchpoints}\}
\end{equation}
and the group
\begin{equation}\label{defsB}
  \sB(d) = \sS(d) \ltimes (\Zt)^d
\end{equation}
is the centralizer of a fixed-point-free involution in $\sS(2d)$. 
This group is also known as the Weyl group of the root systems $B_d$ and $C_d$
and has a convenient realization by signed permutations, 
that is, by automorphisms of $\{\pm 1, \dots, \pm d\}$
commuting with the involution $k\mapsto - k$.

\subsubsection{}\label{sSe}

The group $\sB(d)$ has  a remarkable central 
extension 
\begin{equation}
0 \to \Zt \to \sC(d) \to \sB(d) \to 0
\label{exsec}
\end{equation}
which we will call the \emph{Sergeev 
group} in honor
of A.~Sergeev who recognized in \cite{S1} its importance in 
the representation theory of superalgebras and its
relation to the projective representations of the 
symmetric groups. 

The Sergeev group $\sC(d)$ is the semidirect 
product 
\begin{equation}
  \label{serg}
  \sC(d) \cong  \sS(d) \ltimes \Cliff(d)\,,
\end{equation}
where $\Cliff(d)$ is the Clifford group 
generated by the involutions $\xi_1,\dots,\xi_d$
and and a central involution 
$\epsilon$ subject to the relation
\begin{equation}
  \label{Cliff}
  \xi_i\, \xi_j = \epsilon \,\xi_j \,\xi_i\,, \quad i\ne j \,. 
\end{equation}
The group $\sS(d)$ acts on $\Cliff(d)$ by permuting the 
$\xi_i$'s. Setting $\epsilon=1$, gives the 
group $\sB(d)$. 
\footnote{
Traditionally, the Sergeev group is defined slightly
differently, namely, one sets $\xi_i^2 = \epsilon$ 
instead of $\xi_i^2=1$. While this leads to a nonisomorphic
group, the representation theory is the same.}

\subsubsection{}

The symmetric group $\sS(d)$ is naturally a subgroup of 
$\sB(d)$ and $\sC(d)$. We will call the corresponding
elements of $\sB(d)$ and $\sC(d)$ pure permutations. Let
$$
\delta \in \pi_1(\cBp)
$$
be a loop encircling one of the branchpoints.
By construction, the monodromy $\Mon_{f_\gamma}(\delta)$ is conjugate
in $\sB(d)$ to a pure permutation with odd cycles. 
We have  the following well-known
\footnote{The proof follows from observing that
the commutator of a $k$-cycle $(123\dots k)$ with 
$\xi_1 \xi_2 \cdots \xi_k$ equals $\epsilon^{k-1}$.}

\begin{Lemma}\label{lcan}
Let $g\in\sB(d)$ be conjugate to a pure permutation with 
odd cycles. Then only one of the two preimages of $g$ in $\sC(d)$ 
is conjugate to a pure permutation. 
\end{Lemma}

We will call this distinguished lift of $g\in\sB(d)$ to the Sergeev group 
the \emph{canonical lift} of $g$. 
Given a homomorphism \eqref{phig}, we say that is 
has a canonical lift to the Sergeev group $\sC(d)$,
if there exists a homomorphism 
\begin{equation}
  \label{liftphig}
 \hMon_{f_\gamma} : \pi_1(\cBp) \to  \sC(d)\,,
\end{equation}
that covers $\Mon_{f_\gamma}$, such that for all loops $\delta$ as above 
$\hMon_{f_\gamma}(\delta)$ is the canonical lift of 
$\Mon_{f_\gamma}(\delta)$. We define
\begin{equation}
  \label{defsf}
  \Lift(f_\gamma) = 
  \begin{cases}
    0\,, & \textup{$\Mon_{f_\gamma}$ has a canonical lift}\,, \\
    1\,, & \textup{otherwise} \,. 
  \end{cases}
\end{equation}

\subsubsection{}

In practice, this means the following. The group 
$\pi_1(\cBp)$ is a 
group with one relation of the form 
\begin{equation}
  \label{rel}
  \prod_{i=1}^g [\alpha_i,\beta_i] \prod \delta_i = 1\,,
\end{equation}
where the loops $\alpha_i$ and $\beta_i$ follow the 
standard cycles of $\cB$ and the 
loops $\delta_i$ encircle the branchpoints of $f$. 
The square brackets in \eqref{rel} denote the
commutator
$$
[\alpha,\beta] = \alpha \beta \alpha^{-1} \beta^{-1} \,.
$$
Note that a commutator always has a canonical lift to 
any central extension. 

Since we insist that all 
elements $\hMon_{f_\gamma}(\delta_i)$ are canonical
lifts, this fixes the image of the left-hand side
of \eqref{rel} under the map $\hMon_{f_\gamma}$ completely.
If follows that 
\begin{equation}
  \label{rel2}
  \prod_{i=1}^g \left[
\hMon_{f_\gamma}(\alpha_i), \hMon_{f_\gamma}(\beta_i)\right] 
\prod_i \hMon_{f_\gamma}(\delta_i)  = \epsilon^{\, \Lift(f_\gamma)} \,. 
\end{equation}

\subsubsection{}

Our group-theoretic description of the invariant \eqref{defDq}
is given in the following 

\begin{Theorem}\label{t1}
We have 
\begin{equation}
  \label{thm1}
  \Delta_f (\gamma) = \Lift(f_\gamma)\,,
\end{equation}
or, in other words, 
\begin{equation*} 
\Big( q_{f^*L}(\gamma) - q_L(f_* \gamma) = 0 \Big) \quad 
\Leftrightarrow  \quad 
\left(
\begin{gathered}
 \textup{monodromy of $f_\gamma$ has a } \\
  \textup{canonical lift to Sergeev group} 
\end{gathered}
\right) \,. 
\end{equation*}
\end{Theorem}

\subsection{Torus coverings}

\subsubsection{}

We now specialize the case \eqref{ftoT} of a 
branched covering of the torus. Let the 
ramification data of $f$ be specified by
a collection 
\begin{equation}
  \label{defMu}
  M=(\mu^{(1)},\mu^{(2)},\dots)
\end{equation}
of partitions of $d$ into 
odd parts.
We take $L=0$ and set
$$
p(f)=p(f^*0) \,.
$$
There is a classical method to enumerate
coverings \eqref{ftoT} in terms of 
characters of the symmetric group, 
see for example \cite{J}. Our goal
is to refine this enumeration 
taking into account the parity $p(f)$. 
This will involve characters of 
the Sergeev group. There exists
a very close relation between 
these characters and characters
of projective representations of 
the symmetric group, see \cite{Joz,S1,S2,St,Ya}. 

\subsubsection{}

Introduce the following $\Zt$-grading in
the group $\sC(d)$
\begin{equation}
  \label{grad}
  \deg \xi_i = 1\,, \quad \deg g = \deg \epsilon = 0\,, \quad g\in \sS(d) \,. 
\end{equation}
We will be interested in the 
irreducible $\Zt$-graded modules (also known as 
\emph{supermodules}) of $\sC(d)$ in which the 
central element $\epsilon$ acts nontrivially. 
Such modules are indexed by \emph{strict}
partitions $\lambda$ of $d$. By definition, 
a partition $\lambda$ is strict if
$$
\lambda_1 > \lambda_2 > \dots \ge 0  \,, 
$$
that is, if no nonzero part of $\lambda$ is
repeated. 

Let $V^\lambda$ denote the irreducible 
$\sC(d)$-supermodule corresponding 
to a strict partition $\lambda$. 
Let 
$$
g\in \sS(d)\subset \sC(d)
$$
be a permutation with odd cycles. 
The conjugacy class of $g$ acts in 
$V^\lambda$ as multiplication by 
a constant, known as the 
\emph{central character} of $V^\lambda$. 
It will be denoted by $\fb_\mu(\lambda)$, where $\mu$
is the cycle type of $g$.  
\footnote{
To simplify the
notation, we do not decorate $\fb_\mu(\lambda)$
with any additional symbols to distinguish
it from the identical notation used  for the 
central characters of the symmetric 
group in  e.g.\ \cite{OP}. Same remark 
applies to the functions $\pb_k$ below 
and to algebra $\Lambda$ they generate. Since the symmetric
group analogs of these objects will not
appear in this paper, this should not
lead to confusion.}

\subsubsection{}

It is very convenient to omit parts equal to $1$ 
from the odd partitions $\mu$ in \eqref{defMu} and
in the label of $\fb_\mu(\lambda)$. We also 
set 
$$
\fb_\mu(\lambda)=0\,,\quad |\mu| > |\lambda|\,.
$$
With these conventions, it can be shown that 
$\fb_\mu(\lambda)$ is in fact a \emph{polynomial}
in $\lambda$ in the following precise sense. 
By definition, set
\begin{equation}
  \label{defpb}
  \pb_k(\lambda) = \sum_i \lambda_i^k - \frac{\zeta(-k)}{2} \,,
\quad k=1,3,5,\dots 
\end{equation}
where the constant is introduced for later
convenience (compare with \eqref{Eis}). 
Denote by $\Lambda$ the algebra 
\begin{equation}
  \label{LS}
  \Lambda = \Q[\pb_1,\pb_3,\pb_5,\dots]\,,
\end{equation}
that the polynomials \eqref{defpb} generate. It is 
known as the \emph{algebra of supersymmetric 
functions}.

It is 
known, see Section 6 in \cite{I}, that
$$
\fb_\mu \in \Lambda\,,
$$
and moreover
\begin{equation}
    \label{lfmu}
    \fb_\mu = \frac{1}{z_\mu} \, \pb_\mu + \dots \,,
\end{equation}
Here dots denote lower degree terms, 
$\pb_\mu = \prod \pb_{\mu_i}$, and 
$$
z_\mu = \prod_{k} k^{\rho_k} \, \rho_k !\,,
$$ 
where $\mu=(3^{\rho_3},5^{\rho_5}, \dots)$. 
The formula \eqref{lfmu} is the Sergeev group 
analog of a result of Kerov and Olshanski \cite{KO}.

\subsubsection{}

Our next result gives an enumeration of
nonisomorphic branched covering \eqref{ftoT} with 
given ramification $M$ with 
respect to their parity $p(f)$. Each covering 
appears with a weight which is the 
reciprocal of the order of its 
automorphism group. Such weighting
is standard when counting some 
objects up to isomorphism. 

\begin{Theorem}\label{t2} We have
\begin{equation}
  \label{thm2}
  \sum_{f} \frac{(-1)^{p(f)}\, q^{\deg f}}
{|\Aut(f)|}  = 2^{\chi(\cC)/2} 
\sum_{\textup{strict $\lambda$}} (-1)^{\ell(\lambda)} \, q^{|\lambda|} 
\prod \fb_{\mu^{(i)}}(\lambda) \,, 
\end{equation}
where the summation is over all nonisomorphic 
coverings $f$ with specified ramification $\mu^{(i)}$ 
over fixed points $q_i\in \T^2$. 
\end{Theorem}

Here $\ell(\lambda)$ denotes the length of $\lambda$, that is,
the number of its nonzero parts and $\chi(\cC)$ is the Euler
characteristic of the cover, which can be computed from 
the ramification data \eqref{defMu} as follows 
$$
\chi(\cC) =  \sum_{k} \left(\ell\left(\mu^{(k)}\right)-
\left|\mu^{(k)}\right|\right) \,.
$$

\subsubsection{}

In particular, in the unramified case $M=\emptyset$ we get
the following series 
$$
(q)_\infty = \prod_{n>0} (1-q^n) = 
\sum_{\textup{strict $\lambda$}} (-1)^{\ell(\lambda)} \, q^{|\lambda|} \,.
$$
Set, by definition, 
$$
F_M(q) = \frac1{(q)_\infty} \sum_{f} \frac{(-1)^{p(f)}\, q^{\deg f}}
{|\Aut(f)|} \,, 
$$
where the sum over $f$ is the same as in Theorem \ref{t2}. 
 The series $F_M(q)$ enumerates coverings 
without unramified connected components. 

As a corollary to Theorem \ref{t2}, we obtain the 
following quasimodular property of the series $F_M(q)$. 
Let
\begin{equation}
  \label{Eis}
  E_k(q) = \frac{\zeta(-k+1)}{2} + 
\sum_n \left(\sum_{d |n} d^{k-1}\right) \, q^n\,,
\quad k=2,4,6,\dots\,,
\end{equation}
denote the Eisenstein series of weight $k$. These
series generate the  algebra 
$$
QM_*=\C[E_2, E_4, E_6] 
$$
known as the algebra of \emph{quasimodular forms}. This 
algebra is graded by the weight. We will denote by $QM_k$ the 
graded component of weight $k$ and by $QM_{\le k}$ the 
filtered component of weight at most $k$. 

Define the weight of a partition $\mu$ by
$$
\wt \mu  = |\mu| + \ell(\mu)\,,
$$
where $\ell(\mu)$ is the length of $\mu$, 
and set $\wt M = \sum \wt \mu^{(k)}$. We have the 
following

\begin{Corollary} \label{c1} 
  \begin{equation}
    \label{cor1}
     F_M(q) \in QM_{\le \wt M} \,.
  \end{equation}
\end{Corollary}

The analogous generating series without the 
$(-1)^{p(f)}$ factor is known to be 
quasimodular form of the same weight 
by the result of \cite{BO}. It follows
that the generating series for 
even and odd coverings separately
are quasimodular. 

\subsection{Acknowledgments}

A.~E.\ was partially supported by NSF grant DMS-9704845,
and fellowships from the Sloan and Packard foundations.
A.~O.\ was partially supported by
DMS-0096246 and fellowships from the Sloan and Packard foundations.
R.~P.\ was partially supported by DMS-0071473
and fellowships from the Sloan and Packard foundations.

\section{Proof of Theorem \ref{t1}}

\subsection{}

\begin{figure}[!hbtp]
  \begin{center}
    \scalebox{1}{\includegraphics{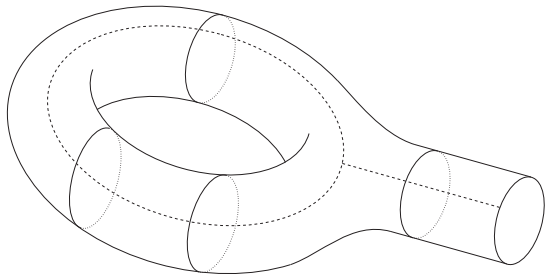}} 
    \caption{A decomposition into tubes and pairs of pants.}
    \label{f1}
  \end{center}
\end{figure}

The proof of Theorem \ref{t1} will be based on 
the usual TQFT technique of cutting the punctured
curve $\cBp$ into 
tubes and pairs of pants,
see for 
example Figure \ref{f1} in which  a punctured torus is
decomposed into 3 tubes and one pair of pants. 
The dashed line in Figure \ref{f1} cuts each 
tube and each pair of pants into a 
simply-connected piece. This dashed line will be 
useful for keeping track of how the
pieces of $\cBp$ fit together and
for understanding the coverings
of $\cBp$. 

By a choosing a suitable 
decomposition of $\cBp$, we will achieve that the
curve $f_*(\gamma)$ takes a particularly 
simple standard form on each piece. 
These standard pieces will have a simple 
interpretation in terms of the Sergeev group $\sC(d)$.

\subsection{}

\begin{figure}[!hbtp]
  \begin{center}
    \scalebox{0.64}{\includegraphics{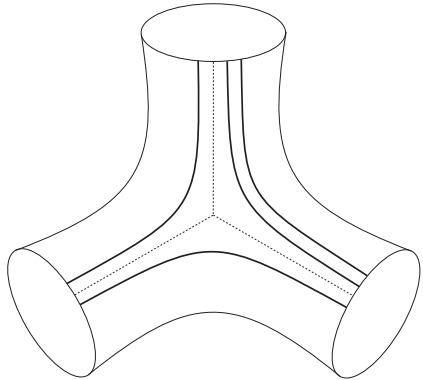}} 
    \caption{The standard pair of pants.}
    \label{f2}
  \end{center}
\end{figure}

\begin{figure}[hbtp]\psset{unit=0.5 cm}
  \begin{center}
    \begin{pspicture}(-11,-3.8)(11,3.8)
    \rput(0,0){\scalebox{0.64}{\includegraphics{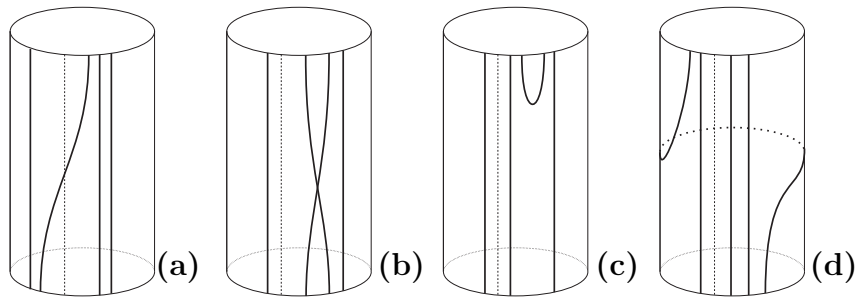}}} 
    \rput[lb](-6.7,-3.5){\textbf{(a)}} 
    \rput[lb](-0.8,-3.5){\textbf{(b)}} 
    \rput[lb](5,-3.5){\textbf{(c)}}
    \rput[lb](10.7,-3.5){\textbf{(d)}}
    \end{pspicture}
    \caption{The standard tubes.}
    \label{f3}
  \end{center}
\end{figure}

Concretely, we require that on all pairs of 
pants the curve $f_*(\gamma)$ has the form 
shown in Figure \ref{f2}, while on the tubes
it has one of the 4 forms shown in Figure \ref{f3}.
Note that, with the exception of case (d) in Figure \ref{f3},
the curve $f_*(\gamma)$ stays in a small neighborhood of 
the dashed line. As a subcase, case (d)
includes a circle around the tube with no vertical 
lines. Also note that the punctures
of $\cBp$ should be capped off by tubes
 of the form (c) with no vertical 
lines remaining.

It is clear from this description that the 
number $\Delta_f (\gamma)$ of new self-intersections
of $f_*(\gamma)$ equals the number of tubes
of type (b) such that the two intersecting 
branches of $f_*(\gamma)$ come from different
sheets of the covering. 
%
%
%
%

\subsection{}

\begin{figure}[hbtp]\psset{unit=0.5 cm}
  \begin{center}
    \begin{pspicture}(0,0)(10,5)
    \rput[lb](0,0){\scalebox{1}{\includegraphics{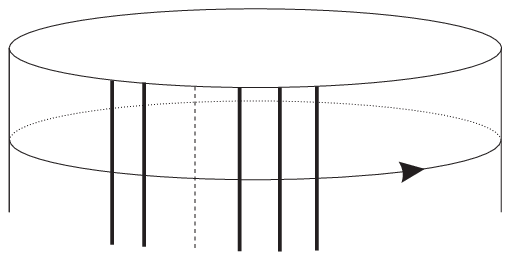}}} 
    \rput[b](2.1,3.6){$i_1$}
    \rput[b](2.9,3.5){$i_2$}
    \rput[b](6.4,3.5){$j_1$}
    \rput[b](5.5,3.5){$j_2$}
    \rput[b](9,1){$\nu$}
    \end{pspicture}
    \caption{Monodromy around the neck $\nu$.}
    \label{f4}
  \end{center}
\end{figure}

We now translate these pictures into the language of 
branched coverings. Consider the
monodromy around the neck $\nu$ of a tube or of a pair 
of paints, see Figure \ref{f4}. We choose the base point
away from the dashed line and study how the 
different sheets of $f_\gamma$ are glued together
across the dashed line. 
%
%
%

Let $g\in \sS(d)$ be the 
monodromy of $f$ around $\nu$, that is, the permutation 
that describes how the sheets of $f$ are 
identified across the dashed line. 
By construction, the monodromy $\Mon_{f_\gamma}(\nu)\in \sB(d)$
is a signed version of $g$. More precisely, 
a sign appears whenever the lift of the 
monodromy path crosses the curve $\gamma$. 
Let us label the branches of $f_*(\gamma)$
by the sheets of $f$ that they come from,
as in Figure \ref{f4}. We choose 
\begin{equation}
  \label{lMnu}
  \hMon_{f_\gamma}(\nu) = \xi_{j_1} \xi_{j_2} \cdots \, g  \, \cdots \xi_{i_2} \xi_{i_1} 
\end{equation}
as the lift of $\Mon_{f_\gamma}(\nu)$ to the Sergeev group. 

For example, the pair of pants from Figure \ref{f2} reflects
the following multiplication rule 
\begin{multline*}
  \left(\xi_{j_1} \xi_{j_2} \cdots \, g  \, \cdots \xi_{i_2} \xi_{i_1}\right)
\cdot \left(\xi_{i_1} \xi_{i_2} \cdots \, h  \, \cdots \xi_{k_2} \xi_{k_1}\right) 
=
\xi_{j_1} \xi_{j_2} \cdots \, g h  \, \cdots \xi_{k_2} \xi_{k_1}\,.
\end{multline*}
Similarly, the type (a) tube corresponds to the identity
$$
g \, \xi_i = \xi_{g(i)} \, g \,,
$$
while the type (c) tube is a geometric version of $\xi_i^2=1$.

\subsection{}

\begin{figure}[hbtp]\psset{unit=0.5 cm}
  \begin{center}
    \begin{pspicture}(0,-0.5)(4,7.8)
    \rput[lb](0,0){\scalebox{0.64}{\includegraphics{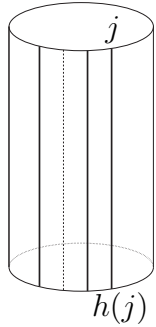}}} 
    \rput(2.8,7){$j$}
    \rput[t](3,0){$h(j)$}
    \end{pspicture}
    \caption{The glue tube.}
    \label{f5}
  \end{center}
\end{figure}

Two necks can be glued together when the two 
monodromies are conjugate. The conjugating element 
gives the gluing rule. This conjugating element 
can be viewed as a new kind of tube. If it is 
a pure permutation $h\in \sS(d)$, then this 
tube looks rather dull, see Figure \ref{f5}, 
because of the relation 
$$
h \cdot \left(\xi_{j_1} \xi_{j_2} \cdots \, g  \, \cdots \xi_{i_2} \xi_{i_1}\right) 
\cdot
h^{-1} = 
\xi_{h(j_1)} \xi_{h(j_2)} \cdots \, h g h^{-1}  \, \cdots \xi_{h(i_2)} \xi_{h(i_1)} \,.
$$
If, however, this conjugating element is $\xi_i$ for some $i$
than the corresponding tube is of type (d). 
%
%

\subsection{}

Finally, the tube of type (b) is a plain violation of the
Sergeev group law, unless the two intersecting branches
come from the same sheet of the covering. Each instance 
of such violation contributes a factor of $\epsilon$ to 
the product \eqref{rel2}. Hence, in total, the product
\eqref{rel2} equals $\epsilon^{\Delta_f(\gamma)}$. 
This completes the proof of Theorem \ref{t1}.

\section{Torus coverings}

\subsection{Proof of Theorem \ref{t2}}

\subsubsection{}

We will begin with a quick review 
of the relation between branched
coverings and characters of groups,
see, for example, \cite{J}.

Let $\sG$ be a finite group and let a conjugacy  class
$C_i\subset \sG$ be specified for each 
element of a finite set of points $\{r_i\}\subset \T^2$. 
A natural question associated to these data is 
to count the homomorphisms
$$
\psi : \pi_1(\T^2\setminus \{r_i\}) \to \sG
$$
sending the conjugacy class of a small loop around
$\{r_i\}$ into the conjugacy class $C_i$.  By definition, set
\begin{equation}
h_\sG(C_1,C_2,\dots) = \frac{\textup{\# of such $\psi$}}{|\sG|}\,.
 \label{h_G}
\end{equation}
In the numerator here, one counts all homomorphisms $\psi$
without any equivalence relation imposed. Alternatively,
the number $h_\sG$ can be interpreted as the 
automorphism weighted count of homomorphisms up the 
action of $\sG$ by conjugation. 

A classical formula, which goes back to 
Burnside, for the number
$h_\sG$ is the following 
\begin{equation}
h_\sG(C_1,C_2,\dots)  = \sum_{\lambda\in \sG^\wedge} \prod \fb_{C_i}(\lambda)\,,
\label{hchar}
\end{equation}
where $\sG^\wedge$ is the set of irreducible complex
representations of $\sG$ and $\fb_C(\lambda)$ is the
central character of the representation $\lambda$.

\subsubsection{} 
The cases of interest to us will be 
\begin{equation}
\sG = \sS(d),\, \sB(d),\, \sB_0(d),\, \sC(d),\, \sC_0(d)\,, 
\label{grps}
\end{equation}
where 
$$
\sC_0(d) \subset \sC(d)\,, \quad \sB_0(d) \subset \sB(d)
$$
are subgroups of 
index two formed by even elements with 
respect to the grading \eqref{grad}. 

The conjugacy classes of interests to us
all come from $\sS(d)$ and are indexed 
by partitions of $\mu$ of $d$ into odd 
parts. We will denote by
$$
h_\sG(M) = h_\sG(C_{\mu^{(1)}},C_{\mu^{(2)}},\dots)
$$
the count of homomorphisms corresponding to
the ramification data \eqref{defMu}.

\subsubsection{}

All groups \eqref{grps} have a natural map to $\sS(d)$ 
and we will denote by 
$$
h_\sG(M;\phi)
$$
the count of those $\psi$ that composed with 
$\sG\to \sS(d)$ yield some fixed homomorphism 
\begin{equation}
  \label{homphi}
  \phi: \pi_1\left(\T^2\setminus\{r_i\}\right) \to \sS(d) \,.
\end{equation}

\begin{Lemma}\label{l1}
For any $\phi$, we have 
\begin{equation}
h_{\sB(d)}(M;\phi) = \frac{1}{d!} \prod_i 2^{2g(\cC_i)-1}\,,\label{hB1}
\end{equation}
where the product ranges over all connected components $\cC_i$ of 
the curve $\cC$ and $g(\cC_i)$ denotes the genus of $\cC_i$. 
\end{Lemma}

\begin{proof}
If $\cC$ is connected, then the set $J_2(\cC)$
is in a natural bijection with the orbits of the
action of $(\Zt)^d \subset \sB(d)$ by conjugation
on all lifts of $\phi$ to $\sB(d)$. Hence, there are
$2^{2g(\cC)}$ such orbits. 

The stabilizer of the $(\Zt)^d$--action is 
isomorphic to $\Zt$ for any orbit. Indeed, any element
in the stabilizer commutes with the subgroup of $\sS(d)$
generated by $\phi$. Since this group is transitive,
it follows that the stabilizer is the group $\Zt$
generated by the involution $k\mapsto -k$. Hence
the number of elements in each orbit is $2^{d-1}$.

Multiplying $2^{d-1}$ by $2^{2g(\cC)}$ and dividing by 
the order of the group gives \eqref{hB1}. 
The generalization to a disconnected curve $\cC$
is straightforward. 
\end{proof}

\subsubsection{}

Our goal now is to specialize Theorem \ref{t1}
to the case $\cB=\T^2$ and $L=0$. 
For $L=0$, the quadratic form $q_L(\nu)$,
$\nu\in H_1(\T^2,\Zt)$,   takes
a very simple form, namely 
\begin{equation}
  \label{q_0}
  q_0(\nu) = 
  \begin{cases}
    0\,, &\nu=0\,, \\
    1\,, &\nu\ne 0\,. 
  \end{cases}
\end{equation}

We now make the following simple observation: 

\begin{Lemma} For $\gamma\in H_1(\cC,\Zt)$, 
the push-forward $f_*\gamma\in H_1(\T^2,\Zt)$  vanishes 
if and only if the monodromy of $f_\gamma$
lies in the group $\sB_0(d)$. 
\end{Lemma}

We now see that Theorem \ref{t1} specializes to 
the 
group-theoretic description 
of the quadratic form $q_{f^*L}$
summarized in the following table:

\renewcommand{\arraystretch}{1.5}
\begin{center}
  \begin{tabular}{|c|c|c|}
\hline 
\multicolumn{2}{|r|}
{Monodromy of $f_\gamma$ lifts to}& $q_{f^*L}(\gamma)$\\
\hspace{0.4cm}  $\sB_0(d)$ \hspace{0.4cm}  & 
$\sC(d)$ &  \\
\hline
yes & yes & 0 \\
yes & no & 1 \\
no & yes & 1 \\
no & no & 0 \\
\hline
  \end{tabular}
\end{center}

\subsubsection{}

This table translates into the following: 

\begin{Proposition}
Let $\phi$
be a homomorphism \eqref{homphi} corresponding
to a covering $f$ as in \eqref{ftoT}. Then
 \begin{multline}
   \label{parphi}
   \frac{(-1)^{p(f)}}{d!}  = {2^{\chi(\cC)/2}}
\Big(h_{\sB(d)}(M;\phi) - h_{\sB_0(d)}(M;\phi) \\ -
h_{\sC(d)}(M;\phi) + h_{\sC_0(d)}(M;\phi)\Big)\,. 
 \end{multline}
\end{Proposition}

\begin{proof}
  Let $\psi$ be a lift of $\phi$ to $\sB(d)$ corresponding
to some class $\gamma\in H_1(\cC,\Zt)$. If the image 
of $\psi$ lies in $\sB_0(d)$, then it contributes twice
as much to $h_{\sB_0(d)}(M;\phi)$ because the order 
of the group is smaller by a factor of $2$. 

Suppose that $\psi$ has a canonical lift to $\sC(d)$. 
Then, in fact, there are $4$ such lifts 
because the monodromy around the two periods
of the torus can be lifted arbitrarily. Since 
$|\sC(d)|\big/|\sB(d)|=2$,
the homomorphism $\psi$ contributes twice as much 
to $h_{\sC(d)}(M;\phi)$ as it does to $h_{\sB(d)}(M;\phi)$. 

Finally, if $\psi$ lifts to $\sC_0(d)$, its contribution 
to $h_{\sC_0(d)}(M;\phi)$ is 4 times as much by the 
combination of two arguments. 

{}Using the above table, one checks that in all cases 
these contribution of $\psi$ 
to the alternating expression in the right-hand
side of \eqref{parphi} combine to  $(-1)^{q_{f^*L}(\gamma)}$. 
Now it remains to apply Lemma \ref{l1} and 
formula \eqref{arf}. 
\end{proof}

\subsubsection{}

It is clear from \eqref{hchar} that
\begin{equation}
  \label{pCmB}
  h_{\sC(d)}(M) - h_{\sB(d)}(M) = 
\sum_{\lambda\in \sC(d)^\wedge_-} \prod \fb_{\mu^{(i)}}(\lambda)\,,
\end{equation}
where $\fb_\mu$ is the central character corresponding to 
an odd partition $\mu$ of $d$ and $\sC(d)^\wedge_-$ denotes
the set of irreducible representations of $\sC(d)$ in which 
the central element $\epsilon$ acts as multiplication by $-1$. 

It follows from \eqref{parphi} and  \eqref{pCmB} that 
\begin{multline}
  \label{enumer1}
  2^{-\chi(\cC)/2} \sum_{\deg f = d } 
\frac{(-1)^{p(f)}}{|\Aut f|} = \\
\left(\sum_{\lambda\in \sC_0(d)^\wedge_-} \prod \fb_{\mu^{(i)}}(\lambda)
\right) - 
\left(\sum_{\lambda\in \sC(d)^\wedge_-} \prod \fb_{\mu^{(i)}}(\lambda)
\right)\,,
\end{multline}
where the summation on the left is over nonisomorphic degree
coverings $f$ with ramification data $M$. 

\subsubsection{}

Consider the algebra 
$$
\sA(d) = \C[\sC(d)]\big /(\epsilon + 1) \,.
$$
The grading \eqref{grad} makes it a
semisimple associative superalgebra. 
The structure of this algebra, or
equivalently, representations of $\sC(d)$
in which $\epsilon$ acts nontrivially, 
were first studied by Sergeev in \cite{S1}
and in many papers since, see the References. 
In particular, the the irreducible 
supermodules $V^\lambda$ of $\sA(d)$ are 
indexed by strict partitions $\lambda$ of $d$.

\subsubsection{}

If $\ell(\lambda)$ is even, then 
the corresponding summand of $\sA(d)$ is isomorphic, 
for some $n$, 
to the associative superalgebra 
\begin{equation}
  \label{gl}
  \gl(n,n) = \left\{
    \begin{pmatrix}
      A_{00} & A_{01} \\
      A_{10} & A_{11}
    \end{pmatrix}
\right\}\,, \quad \deg A_{ij}= i+j\,,
\end{equation}
of all endomorphisms of the 
$\Zt$-graded vector space $V^\lambda$. In particular, 
$V^\lambda$ is an irreducible representation of $\sC(d)$,
which splits into two non-equivalent 
irreducible modules 
$$
V^\lambda=V^\lambda_0 \oplus V^\lambda_1
$$ 
as a module over the even subalgebra
$$
\sA \supset \sA_0 = \C[\sC_0(d)]\big /(\epsilon+1) \,.
$$
However, these two representations, being 
conjugate under the action of $\sC(d)$, have
the same central character when restricted
to permutations with odd cycles. Indeed, any permutation with odd cycles
commutes with some odd element of $\sC(d)$, such as, 
for example, $(123)$ commutes with $\xi_1 \xi_2 
\xi_3$. 

It follows that when $\ell(\lambda)$
is even the contribution of $V^\lambda$
to the first sum in \eqref{enumer1} is twice its 
contribution to the second sum in \eqref{enumer1}.

\subsubsection{}

If $\ell(\lambda)$ is odd then the corresponding 
summand of $\sA(d)$ is  is isomorphic, 
for some $n$, to the associative superalgebra 
\begin{equation}
  \label{fq}
  \fq(n) = \left\{
    \begin{pmatrix}
      A_{0} & A_{1} \\
      A_{1} & A_{0}
    \end{pmatrix}
\right\}\,, \quad \deg A_{i}= i\,. 
\end{equation}
It follows that the corresponding supermodule 
$V^\lambda$ is a direct sum of two non-equivalent 
representations of $\sC(d)$, which however,
become equivalent after restriction to 
$\sC_0(d)$. Therefore, the contribution of
such $\lambda$ to the second sum  
in \eqref{enumer1} is precisely twice its 
contribution to the first sum in \eqref{enumer1}.

\subsubsection{}

We can summarize this discussion as follows
\begin{equation}
  \label{enumer2}
   \sum_{\deg f = d} 
\frac{(-1)^{p(f)}}{|\Aut f|} = 
2^{\chi(\cC)/2} \sum_{\textup{strict partitions $\lambda$ of $d$}} 
(-1)^{\ell(\lambda)} \, 
\prod \fb_{\mu^{(i)}}(\lambda)\,,
\end{equation}
which is equivalent to the statement of 
Theorem \ref{t2}.

\subsection{Quasimodularity of $F_M(q)$}

\subsubsection{}

For $f\in \Lambda$, let
$$
\lan f \ran = \frac1{(q)_\infty} \sum_{\textup{strict $\lambda$}} 
(-1)^{\ell(\lambda)} \, q^{|\lambda|} \, f(\lambda) 
$$
denote the average of $f$ with the respect to the weight 
$(-1)^{\ell(\lambda)} \, q^{|\lambda|}$ on strict 
partitions. To prove Corollary \ref{c1}, it 
clearly suffices to show that 
\begin{equation}
  \label{avp}
  \lan \pb_\mu \ran \in QM_{\wt \mu} \,.
\end{equation}
The proof of \eqref{avp} is straightforward and 
follows \cite{BO}.

\subsubsection{}

Let $t_1,t_3,t_5,\dots$ be formal variables. We have
\begin{multline}
  \label{idenq}
 \frac1{(q)_\infty}
 \sum_{\textup{strict $\lambda$}} \, (-1)^{\ell(\lambda)} \, q^{|\lambda|} \,
 \exp \left(\sum t_k \, \pb_k(\lambda)\right) = \\
\exp\left(-\frac12 \sum t_k \, \zeta(-k)\right) 
\prod_{n>0} \frac{\left(1- q^n \exp\left(\sum t_k \, n^k\right)\right)}
{\left(1- q^n \right)}\,.
\end{multline}
This is a generating function for the averages \eqref{avp}. 
The logarithm of \eqref{idenq} equals
\begin{equation}
  \label{log}
-\frac12 \sum t_k \, \zeta(-k)
  - \sum_{m,n} \frac{q^{mn}}{m}
\left[\exp\left(m \sum t_k \, n^k\right) - 1 \right] 
 \,.
\end{equation}
For an odd partition $\mu$, 
the coefficient of $\prod  t_{\mu_i}$, in the 
expansion of \eqref{log} equals, up to a 
constant factor, 
\begin{equation*}
  \delta_{\ell(\mu),1} \, \frac{\zeta(-\mu_1)}{2}+
 \sum_{n,m} q^{mn} m^{\ell(\mu)-1} n^{|\mu|} 
= 
\left(q \, \frac{d}{dq}\right)^{\ell(\mu)-1} E_{|\mu|-\ell(\mu)+2} 
\in
QM_{\wt \mu}\,. 
\end{equation*}
This is because the operator $q \, \frac{d}{dq}$ preserves
the algebra $QM_*$ and increases weight by $2$, which is 
proved in \cite{KaZ} and can be also seen, for example, from the
heat equation for the genus 1 theta-functions.  

Corollary \ref{c1} now follows by exponentiation.

\end{document}